\documentclass[12pt, reqno]{amsart}
\usepackage{latexsym, amsmath,amsthm,amssymb, amsfonts, graphics}

\theoremstyle{plain}
\newtheorem{thm}{Theorem}[section]
\newtheorem{prop}[thm]{Proposition}
\newtheorem{lem}[thm]{Lemma}
\newtheorem{cor}[thm]{Corollary}

\theoremstyle{definition}
\newtheorem{rem}[thm]{Remark}
\newtheorem{defn}[thm]{Definition}
\newtheorem{eg}[thm]{Example}
\numberwithin{equation}{section}

\newcommand{\C}{\mathbb{C}}
\newcommand{\R}{\mathbb{R}}

\newcommand{\cg}{\mathcal{G}}

\newcommand{\co}{\mathcal{O}}

\newcommand{\T}{\mathrm{T}}
\newcommand{\ud}{\mathrm{ \, d}}

\newcommand{\III}{\mathrm{ I \! I \! I \, } }

\DeclareMathOperator{\GL}{\mathrm{GL}}
\DeclareMathOperator{\SL}{\mathrm{SL}}

\DeclareMathOperator{\SO}{\mathrm{SO}}

\DeclareMathOperator{\sli}{\mathrm{sl}}

\DeclareMathOperator{\Tr}{\mathrm{Tr}}
\DeclareMathOperator{\Res}{\mathrm{Res}}
\DeclareMathOperator{\Ima}{\mathrm{Image}}
\DeclareMathOperator{\Ker}{\mathrm{Kernel}}
\DeclareMathOperator{\diag}{\mathrm{diag}}

\begin{document}

\title{Tzitz\'eica transformation is a dressing action}

\author{\bigskip Erxiao Wang}
\address{University of Texas, Austin, TX 78712, USA}
\email{ewang@math.utexas.edu}

\bigskip
\begin{abstract}
\bigskip

We classify the simplest rational elements in a twisted loop group,
and prove that dressing actions of them on proper indefinite affine
spheres give the classical Tzitz\'eica transformation and its dual.
We also give the group point of view of the Permutability Theorem,
construct complex Tzitz\'eica transformations, and discuss the group
structure for these transformations.
\end{abstract}

\maketitle

\clearpage

\section{Introduction}

In 1910, Tzitz\'eica published a classical paper \cite{Tzi10} on
hyperbolic surfaces in $\R^{3}$ whose Gauss curvature at any point
$p$ is proportional to the fourth power of the distance from a fixed
point to the tangent plane at $p$. He proved
\begin{equation} \label{eqtz}
w_{xy}=e^{w}-e^{-2w}
\end{equation}
is the structure equation, and also constructed a geometric
transformation of such surfaces that is similar to the well-known
B\"acklund transformation of surfaces with constant negative
curvature. These surfaces are invariant under affine
transformations, and they are now known as (proper) affine spheres
in affine differential geometry.

The classical Tzitz\'eica equation \eqref{eqtz} was rediscovered in
many mathematical and physical contexts afterwards (see, e.g.,
\cite{BulDod77}, \cite{Dun02}, \cite{Gaf84}). In recent years,
techniques from soliton theory have been applied to this equation
extensively by, e.g.: Rogers \& Schief in the context of gas
dynamics (\cite{RogSch94}), Kaptsov \& Shan'ko on multi-soliton
formulas (\cite{KapSha97}), Dorfmeister \& Eitner on Weierstrass
type representation (\cite{DorEit01}), and Bobenko \& Schief on its
discretizations (\cite{BobSch99a}, \cite{BobSch99b}). Terng \&
Uhlenbeck (\cite{TerUhl00a}) gave a systematic method to construct
B\"acklund-type transformations via dressing actions of simple
rational loop group elements. It is natural to ask whether the
classical Tzitz\'eica transformation is a dressing action of some
loop element, whether there are new transformations of affine
spheres, and what is the group structure of these transformations.
This paper answers these questions.

In section $2$, we give a brief review of classical results and
provide the Lax pair of the structure equations. In section $3$, we
review the reality conditions for this Lax pair  and give the loop
group description of indefinite affine spheres. We then classify the
simplest rational elements in this loop group and compute their
dressing actions on affine spheres in section $4$. These rational
elements cannot be constructed by projections as in \cite{TerUhl00a}
and the computation is harder. It turns out that one class of
dressing action provides exactly the Tzitz\'eica transformation and
the other provides the dual transformation. In section $5$, we
present the group point of view of the classical permutability
theorem, construct complex Tzitz\'eica transformations and discuss
the group structure of these transformations. Some examples are
presented in the last section.

\medskip
\section{Indefinite affine sphere and its Lax representation}

Classical affine differential geometry studies the properties of
surfaces in $\R^{3}$ invariant under the (equi-)affine
transformations $x \to Ax+v$, where $A \in \SL(3,\R)$ and $x, ~ v
\in \R^{3}$. There are three fundamental affine invariants: the
affine (or Blaschke) metric, the Fubini-Pick cubic form, and the
third fundamental form (or the affine shape operator). These
invariants satisfy certain compatibility equations and the
Fundamental Theorem states that they then determine a surface
uniquely up to affine transformations. Let us first review the
definitions of these invariants (for more details see, e.g.,
\cite{Bla23}, \cite{NomSas94}, \cite{Ter83}). The reader may also
refer to \cite{Bla23}, \cite{BobSch99a}, \cite{BobSch99b} or
\cite{SimWan93} for an elementary description of affine spheres.

Let $X: M \hookrightarrow \R^{3}$ be an immersed surface with
non-degenerate second fundamental form. Let $E = (e_1,e_2, e_3 ) $
be a local $\SL(3,\R)$-frame on $M$ such that $e_1, e_2$ are tangent
to $M$, and $e_3$ is transversal to $M$. Let $\omega_1,\omega_2$
denote the dual coframe of $e_1, e_2$, i.e.,
\[
\ud X =  e_1 \otimes \omega_1 + e_2 \otimes \omega_2.
\]
Let $(\omega_{AB} )$ denote the $\sli (3,\R)$-valued $1$-form
$E^{-1} \ud E$, i.e.,
\[
\ud e_A = \sum_{B=1}^{3} e_B \otimes \omega_{BA}.
\]
Then we have the structure equation:
\begin{equation} \label{eqci}
\begin{cases}
\ud \omega_A & = -\sum_{B} \omega_{AB} \wedge \omega_B , \cr \ud
\omega_{AB} & = -\sum_{C} \omega_{AC} \wedge \omega_{CB}.
\end{cases}
\end{equation}
Since $\omega_{3} = 0$ on $M$, \eqref{eqci} implies that for
$i=1,2$:
\begin{equation} \label{eqcp}
\omega_{3, i} = h_{i1} \omega_1 + h_{i2} \omega_2, \textrm{  with  }
h_{ij}=h_{ji}.
\end{equation}
A direct computation shows that the quadratic form
\begin{equation} \label{eqcn}
g := |\det(h_{ij})|^{-\frac{1}{4}} \sum_{i,j=1}^{2} h_{ij} \omega_i
\omega_j
\end{equation}
is invariant under change of affine frames, and it is called the
\emph{affine metric} of $M$. $M$ is said to be {\it definite\/} or
{\it indefinite\/} if the affine metric is definite or indefinite
respectively.

The {\it affine normal} is $ \, \xi := \triangle X \, / \, 2 $,
where $\triangle$ is the Laplacian of $g$. It satisfies two natural
geometric conditions:
\begin{itemize}
\item[(i)] $\ud \xi ( \cdot ) \in \T M $,
\item[(ii)] $ i_{\xi} \ud V = dvol_{g} $ (the volume form of $g$);
\end{itemize}
and is essentially determined by them.

Take the exterior differentiation of \eqref{eqcp} to get
\begin{equation} \label{eqcr}
\sum_{j} ( \ud h_{ij} + h_{ij} \omega_{3, 3} - h_{ik} \omega_{kj} -
h_{kj} \omega_{ki} ) \wedge \omega^j=0 ,
\end{equation}
and define $h_{ijk}$ by
\begin{equation} \label{eqcq}
\sum_{k} h_{ijk} \omega_k = \ud h_{ij} + h_{ij} \omega_{3, 3} -
h_{ik} \omega_{kj} - h_{kj} \omega_{ki}.
\end{equation}
Then \eqref{eqcp} and \eqref{eqcr} imply that $h_{ijk}$ is symmetric
in $i,j,k$. {\it The Fubini-Pick cubic form\/} is defined as
\[
J := \sum_{i,j,k} h_{ijk} \omega_i \omega_j \omega_k,
\]
which is an affine invariant.

We choose $e_3 = \xi$. Then $\omega_{3,3} = 0$. Exterior
differentiate it to get
\[
 \omega_{13} \wedge \omega_{31} + \omega_{23} \wedge \omega_{32} = 0.
\]
Thus the following form
\[
\III := |\det(h_{ij})|^{ \frac{1}{4} } ( \omega_{1 3} \omega_{3 1} +
\omega_{2 3} \omega_{3 2} )
\]
is symmetric. This is the \emph{third fundamental form}.
Equivalently, we can first define the \emph{affine shape operator}
$S$:
\[
 S (u) := \ud \xi (u), \quad  \forall \, u \in \T_{p}M ,
 \]
 then $ \III(u, v) = g(S(u), v) = g(u, S(v)) $. The {\it affine mean curvature\/} $H$ and the {\it affine Gauss curvature\/} $K$ are defined as $ H = \Tr S \, / \, 2 , ~  K = \det S $.

\begin{defn}
An \emph{affine sphere} is a surface all of whose affine normals
meet at a common point.
\end{defn}

An equivalent definition is $S = H \cdot \mathrm{Id} $, i.e., the
shape operator is a scalar multiple of the identity map and $H$ is
then the affine mean curvature. It follows from the structure
equations that $H$ must be constant. When $H = 0$, all affine
normals are parallel and the center is at infinity. Such surface is
called \emph{improper affine sphere} and has been completely
classified in \cite{Bla23}. When $H \neq 0$, it is called
\emph{proper affine sphere} and we can move the center to the origin
and normalize $H$ to $1$ by scaling the ambient space and changing
the orientation if necessary. Then $ e_3 = \xi =  X$.

\smallskip
From now on we will only consider proper indefinite affine spheres
in $\R^{3}$ with $\xi = X$. First note that there exists local
asymptotic coordinate system $(x,y)$ and a smooth function $w(x, y)$
such that the affine metric is:
\[
g =  e^{w} ( \ud x \otimes \ud y + \ud y \otimes \ud x) .
\]
We choose a frame $e_1 = X_x$, $e_2 = e^{-w} X_y$, and $e_3 = \xi =
X $. Then $ \det(e_1, e_2, e_3) = 1 $, and
\[
\begin{array}{ll}
 \omega_1 = \ud x = \omega_{13} = \omega_{32},  \quad & \omega_2 = e^{w} \ud y = \omega_{23} = \omega_{31}, \\
 \omega_{33} = 0, \qquad \qquad \qquad & (h_{ij}) = \begin{pmatrix}   0   &  1  \\    1  &  0 \end{pmatrix} .
\end{array}
\]
A direct computation using the formula \eqref{eqcq} shows
\[
\begin{cases}
\omega_{21}  =  a \ud x, ~ \omega_{12} = b e^{-2w} \ud y  \quad & \textrm{ for some functions }a, b; \\
J  =  -2a \ud x^{3} - 2b \ud y^{3} . &
\end{cases}
\]
Finally from $ \ud \omega_A + \sum_{B} \omega_{AB} \wedge \omega_B =
0 $ we get
\[
\omega_{11} = - \omega_{22} = w_{x} \ud x .
\]
We have obtained the flat $\sli(3,\R)$-valued $1$-form
\[
\omega = E^{-1} \ud E = \begin{pmatrix}
w_x \ud x   & b e^{-2w} \ud y   & \ud x             \\
a \ud x         & - w_x \ud x       & e^{w} \ud y   \\
e^{w}\ud y  & \ud x                 & 0
\end{pmatrix} .
\]
The compatibility equations are
\begin{equation} \label{eqgc}
 \ud \omega + \omega \wedge \omega = 0 ~ \iff ~
  \begin{cases}
     w_{xy}  = e^{w} - ab \  e^{-2w}, \\
     a_{y}  = 0, \quad b_{x} = 0 .
  \end{cases}
\end{equation}
When $ab \neq 0$ we may reparametrize the asymptotic coordinates to
make $a=b=1$. Then \eqref{eqgc} is simplified to the classical
Tzitz\'eica equation.

\begin{rem}
It is known that ruled proper indefinite affine spheres correspond
to the case $ ab \equiv 0 $ and they have been well understood (see
\cite{NomSas94}). For non-ruled case, the points at which $ab=0$ are
called \emph{planar points}.
\end{rem}

The following observation is crucial for the integrability of proper
indefinite affine spheres: The system \eqref{eqgc} is invariant
under the transformation
\[
a \longrightarrow \gamma a, ~ \qquad ~ b \longrightarrow \gamma^{-1}
b
\]
with $\gamma \in \C \setminus \{0\}$. Thus a family of flat
connections is obtained:
\[
\omega_{\gamma} =  \begin{pmatrix}
w_x \ud x           & \gamma^{-1} b e^{-2w} \ud y  & \ud x  \\
\gamma a \ud x  & - w_x \ud x                      & e^{w} \ud y  \\
e^{w}\ud y          & \ud x                                & 0
\end{pmatrix} .
\]
The zero curvature equation of $\omega_{\gamma}$ is called the
\emph{Lax representation} of \eqref{eqgc}.

When we solve $E_{\gamma}$ from
\begin{equation} \label{eqln}
E_{\gamma}^{-1} \ud E_{\gamma} = \omega_{\gamma}
\end{equation}
for any $\gamma \in \R \setminus \{ 0 \}$, the last column of
$E_{\gamma}$ gives a family of affine spheres, whose affine
fundamental invariants are:
\[
g = 2 e^{w} \ud x \ud y, \qquad S = \mathrm{Id}, \qquad J = -2
\gamma a \, \ud x^{3} - \frac{2}{\gamma} b \, \ud y^{3}.
\]

Let us recall the classical \emph{duality relation} for indefinite
affine spheres. Let $h = e^{\omega}$. Then Tzitz\'eica equation
becomes
\begin{equation} \label{eqtz2}
(\ln h)_{xy} = h - \frac{1}{h^{2}}, \qquad \textrm{or }~ h_{xy}h -
h_{x}h_{y} = h^{3} - 1 .
\end{equation}
Classically \eqref{eqln} with $a=b=1$ was written as a linear system
for $X$:
\begin{equation} \label{eqmap}
\begin{cases}
  X_{xx} = \frac{h_x}{h} X_{x} + \frac{\gamma}{h} X_{y} ~, \\
  X_{xy} = h X  ~, \\
  X_{yy} = \frac{1}{\gamma h} X_{x} + \frac{h_y}{h} X_{y} ~.
\end{cases}
\end{equation}

If $X$ solves \eqref{eqmap}, then we can check that
\[
X^{\ast} := \frac{1}{h} X_{x} \times X_{y}
\]
is a solution of \eqref{eqmap} with $\gamma$ replaced by $ - \gamma
$. Here $\times$ is the vector cross product in $\R^{3}$. This is
clearly a duality relation: $(X^{\ast})^{\ast} = X $.

\smallskip
Finally let us recall the classical Tzitz\'eica transformation:
\begin{thm}[ \cite{Tzi10}] \label{thmctt}
Given a solution $ (h, X) $ of \eqref{eqtz2} and \eqref{eqmap}, and
$\phi_{1}$ any scalar solution of \eqref{eqmap} with parameter
$\gamma_{1}$, then the following transformation produces a new
solution $ (h_1, X_1) $ of \eqref{eqtz2} and \eqref{eqmap}:
\begin{equation} \label{eqtzt}
 \begin{cases}
  h_{1}  := h - 2 (\ln \phi_{1})_{xy} ~ , \\
  X_{1}  := \frac{\displaystyle  (\gamma - \gamma_{1} ) h X
        - 2 \gamma (\ln\phi_{1})_{x} X_{y}
        + 2 \gamma_{1} (\ln \phi_{1})_{y} X_{x} }
        {\displaystyle  (\gamma + \gamma_{1} ) h } ~ .
 \end{cases}
\end{equation}
\end{thm}

\medskip
\section{The reality conditions and loop group description}

Henceforth we assume $a=b=1$. To further reveal the hidden symmetry,
let $\lambda = \sqrt[3]{\gamma}$ and change the frame $E_{\gamma} $
to
\[
F_{\lambda} = E_{\gamma} \diag (1/ \lambda, \lambda, 1) .
\]
The gauged family of flat connections is then
\begin{equation} \label{eqframe}
\theta_{\lambda} = F_{\lambda}^{-1} \ud F_{\lambda}  =
\begin{pmatrix} w_x& 0 &\lambda \cr \lambda & -w_x & 0 \cr 0 & \lambda &0 \end{pmatrix}   \ud x
 + \lambda^{-1} \begin{pmatrix} 0& e^{-2w}& 0 \cr 0& 0 & e^w \cr e^w &0&0 \end{pmatrix}  \ud y .
\end{equation}

For any $g \in \SL(3,\C)$, we need to define $\tau(g):=\bar g$ and
define $\sigma$ by:
\[
\sigma(g):= T \ (g^t)^{-1} T^{-1}, \quad \textrm{ where   } T =
\begin{pmatrix} 0& 1 & 0 \cr - \epsilon  &0&0 \cr 0 & 0 & \epsilon^2 \end{pmatrix} , \quad \epsilon = e^{ \pi {\rm i} / 3 } .
\]
The automorphism $\sigma$ has order $6$ and induces the following
automorphism (still denoted by $\sigma$) on the Lie algebra
$\sli(3,\C)$: $\sigma(A)= -T  A^t T^{-1}$.  Therefore $\sigma$ gives
the eigenspace decomposition: $\sli(3,\C) = \oplus_{j=0}^{5} \cg_j$,
where $\cg_j$ is of eigenvalue $\epsilon^j$. We compute that $X_j
\in \cg_j$ if and only if
\begin{align*}
& X_0=\begin{pmatrix} x_{11} &0&0\cr 0&-x_{11} &0\cr 0&0&0\cr
\end{pmatrix},\ \
    X_1=\begin{pmatrix} 0&0& x_{13} \cr x_{21} &0&0\cr 0& x_{13} &0\cr \end{pmatrix},   \cr
& X_2=\begin{pmatrix} 0&0&0\cr 0&0& x_{23} \cr -x_{23} &0&0
\end{pmatrix}, \ \
    X_3=\begin{pmatrix} x_{11} &0&0\cr 0& x_{11} &0\cr 0&0&-2 x_{11} \end{pmatrix},    \cr
& X_4=\begin{pmatrix} 0&0& x_{13} \cr 0&0&0\cr 0& -x_{13} &0
\end{pmatrix},  \ \
    X_5=\begin{pmatrix} 0& x_{12} &0\cr 0&0& x_{23} \cr x_{23} &0&0 \end{pmatrix}.   \cr
\end{align*}
Note that $\sigma\tau=\tau^{-1}\sigma^{-1}$ implies that
$\oplus_{j=0}^{5} (\sli(3,\R)\cap \cg_j)$ is the corresponding
eigenspace decomposition of $\sli(3,\R)$.

The $\theta_\lambda$ in \eqref{eqframe} satisfies two reality
conditions (first given in  \cite{Mik81}):
\begin{equation}\label{eqre}
 \tau (\theta_\lambda) = \theta_{ \bar{\lambda} }, \qquad
 \sigma (\theta_\lambda) = \theta_{ \epsilon \lambda }.
\end{equation}
When we solve $F_{\lambda}$ in \eqref{eqframe} uniquely with the
initial condition $F(0, 0, \lambda) = I$, it is easy to show that
$F_{\lambda}$ also satisfies the reality conditions \eqref{eqre}.

Let $\C_{\ast}:=\C\setminus\{0\}$. We adopt the following notations
for loop groups:
\begin{eqnarray*}
\Lambda G &=& \{  \textrm{ holomorphic maps from }
         \C_{\ast} \cap (\co_{r} \cup \co_{1/ r})  \textrm{ to } G \ \} , \\
\Lambda_{+} G &=& \{  \textrm{ holomorphic maps from } \C_{\ast}  \textrm{ to } G \ \} , \\
\Lambda_{-} G &=& \{  \textrm{ holomorphic maps } f  \textrm{ from
}  \co_{r} \cup \co_{1/ r}  \textrm{ to  } G   ~ \textrm{with }
f(\infty) =  I \ \} ,
\end{eqnarray*}
where $0< r < 1$ is sufficiently small and
\begin{equation*}
\co_{r} = \{ \lambda \in \C ~:~  | \lambda | < r \}  , \qquad
\co_{1/r} = \{ \lambda \in \C \cup \{ \infty \} ~:~  | \lambda | >
1/ r \} ~ .
\end{equation*}

Similar notations also apply to their Lie algebras. Let
$\Lambda^{\tau,\sigma}G$ denote the subgroup of $g\in \Lambda G$
satisfying the reality conditions \eqref{eqre}. Then
$\theta_\lambda$ in \eqref{eqframe} is a $\Lambda_{+}^{\tau,\sigma}
\sli(3,\C) $-valued flat connection, and the corresponding frame
$F_\lambda$ for indefinite affine spheres lies in
$\Lambda_{+}^{\tau,\sigma} \SL(3,\C) $. Conversely, given any smooth
map $F$ from a domain in $\R^{2}$ to $\Lambda_{+}^{\tau,\sigma}
\SL(3,\C) $ satisfying
\begin{equation} \label{eqloop}
F^{-1}F_{x} = A \lambda + B, \qquad F^{-1}F_{y} = C \lambda^{-1} + D
\end{equation}
with $A_{32} C_{31} \neq 0 $, the last column of $F$ then gives an
affine sphere with $h=A_{32} C_{31}$ and $F$ differs from
$F_{\lambda}$ in \eqref{eqframe} by a simple  gauge. This is the
loop group description for indefinite affine spheres (for details,
see \cite{BobSch99a}, \cite{DorEit01}).

\medskip
\section{Dressing actions of simple rational elements}

Let us briefly review the method of dressing action (the original
idea went back to \cite{ZakSha79} but see \cite{Gue97} or
\cite{Ter03} for an elementary introduction). Let $G = \SL(3,\C)$, $
g(\lambda) \in \Lambda^{ \tau, \sigma}_{-} G $, and $F(x,y, \lambda)
\in \Lambda^{\tau, \sigma}_{+} G $ the frame of an associated family
of indefinite affine spheres. Assume we can do the following
factorization for each fixed $(x,y)$:
\begin{equation} \label{eqfact}
g(\lambda) ~ F(x,y, \lambda) = \tilde{F} (x, y, \lambda) ~ \tilde{g}
(x, y, \lambda) ~ ,
\end{equation}
with $\tilde{F} \in \Lambda^{\tau, \sigma}_{+} G $ and $ \tilde{g}
\in \Lambda^{ \tau, \sigma}_{-} G $. Then $\tilde{F}$ also satisfies
\eqref{eqloop} and generates new affine spheres. We sketch the proof
here. It suffices to prove that $\tilde{F}^{-1} (\tilde{F})_{x}$ and
$\tilde{F}^{-1} (\tilde{F})_{y}$ are linear in $\lambda$ and
$\lambda^{-1}$ respectively. But
\begin{eqnarray*}
\tilde{F}^{-1} (\tilde{F})_{x} &=&  \tilde{g} F^{-1} g^{-1} (g F \tilde{g}^{-1})_{x} \\
                &=& \tilde{g} (F^{-1}F_{x}) \tilde{g}^{-1} + \tilde{g}(\tilde{g}^{-1})_{x} \\
                &=& \tilde{g} (A \lambda +B) \tilde{g}^{-1} - \tilde{g}_{x} \tilde{g}^{-1} .
\end{eqnarray*}
On the left hand side, $\tilde{F}^{-1} (\tilde{F})_{x}$ is
holomorphic in $\lambda \in \C \setminus \{ 0\}$; on the right it
has a simple pole at $\infty$ since $g(\infty)=\tilde{g}(\infty)=I$.
So
\[
\tilde{F}^{-1} (\tilde{F})_{x} = \widetilde{A} \lambda +\tilde{B} .
\]
Similarly, $\tilde{F}^{-1} (\tilde{F})_{y}$ is linear in
$1/\lambda$. This completes the proof.

Furthermore, $g \ast F := \tilde{F}$ defines a group action of
$\Lambda^{\tau, \sigma}_{-} G $ on the frames of affine spheres,
which is called the \emph{dressing action}.

The factorization \eqref{eqfact} can indeed be done on a dense open
subset of $\Lambda^{\tau, \sigma} G $ (see \cite{BerGue91},
\cite{PreSeg86}). There is no explicit construction for general $g$,
but when $g$ is a rational element, the factorization can be carried
out using residue calculus (see \cite{TerUhl00a}). In search of
simple rational elements in $\Lambda^{ \tau, \sigma}_{-} \SL(3,\C)
$, it helps to write $\sigma=\nu\circ\mu$ as the composition of two
commuting automorphisms, where
\begin{eqnarray*}
 \nu(g) &:=& Q g Q^{-1} \qquad  \text{with }  Q=\diag(\epsilon^{4}, {\epsilon}^2, 1), \\
 \mu(g) &:=& P (g^t)^{-1} P \quad  \text{with } P= \begin{pmatrix} 0&1&0 \\ 1&0&0 \\ 0&0&1 \end{pmatrix}.
\end{eqnarray*}
Here $\nu$ has been called the Coxeter-Killing automorphism, and the
involution $\mu$ is the unique outer automorphism of $\SL(3,\C)$
modulo inner ones. We observe that  an element $g(\lambda) \in
\Lambda G$ lies in $\Lambda^{\tau, \sigma}G $ if and only if
\begin{equation} \label{eqrea}
\tau(g_{\bar{\lambda}} ) = g_\lambda, \qquad  \nu (g_{\lambda}) =
g_{\epsilon^4 \lambda}, \qquad  \mu (g_{\lambda})  = g_{- \lambda} .
\end{equation}

\begin{rem}
$\tau$ and $\mu$ define a symmetric space $\SL(3,\R) / \SO(2,1) $.
Here  $\SO(2,1)$ is the isometry group of the quadratic form given
by the symmetric matrix $P$, i.e., $2 x_1x_2 + x_3^2$ on $\R^3$. It
was proved in \cite{Ter03} that if $F_{\lambda}$ is the  frame of an
associated family of indefinite affine spheres, then $F_{-1}
F_1^{-1}$ is a harmonic map from $\R^{1,1}$ to the symmetric space
$\SL(3,\R) / \SO(2,1) $.
\end{rem}

We will first study rational elements in $\Lambda^{\sigma}_{-}G$.
Due to the $\nu$-reality condition in \eqref{eqrea}, the simplest
rational element in $\Lambda^{\sigma}_{-} G$ may have only $3$
simple poles $\{ \alpha, \epsilon^{2} \alpha, \epsilon^{4} \alpha
\}$ with $\alpha \in \C_{\ast}$. The element can always take the
following special form:
\begin{equation} \label{eqsf}
g(\lambda) = I + \frac{2\alpha }{\lambda - \alpha }A +
\frac{2\epsilon^2 \alpha }{\lambda - \epsilon^2 \alpha }B +
\frac{2\epsilon^4 \alpha }{\lambda - \epsilon^4 \alpha }C .
\end{equation}
Plug it into $(\nu,\mu)$-reality conditions in \eqref{eqrea} and
compare the residues at each pole, we obtain that $g \in
\Lambda^{\sigma}_{-}G$ if and only if
\begin{equation}
\left\{
\begin{array}{l}
 B=Q^{-1}AQ, \qquad C=QAQ^{-1} , \\
 A^{t} P ( I - A - 2 \epsilon B + 2 \epsilon^{2} C ) = 0 .
\end{array} \right.
\end{equation}
Write $A=(a_{ij})$, and we compute that
\begin{equation}\label{abc}
( I - A - 2 \epsilon B + 2 \epsilon^{2} C ) = \begin{pmatrix}
 1-3a_{11} & -3a_{12} & 3a_{13}    \\
 3a_{21} & 1-3a_{22} & -3a_{23}    \\
 -3a_{31} & 3a_{32} & 1-3a_{33}
\end{pmatrix} .
\end{equation}
If the rank of $A$ is $3$, we get $A= I /3$ and $g(\lambda)$ is
trivial. So there are two types left for $A$: rank $1$ type and rank
$2$ type. A long but not hard computation implies that the
\textbf{rank $1$ type} is as follows:
\begin{equation}\label{res1}
A=\frac{1}{3} \begin{pmatrix} \tfrac{b}{2ab-1} \\ a \\ 1
\end{pmatrix}
\begin{pmatrix}a&b&1 \end{pmatrix} = \frac{1}{3} \begin{pmatrix}
 \tfrac{ab}{2ab-1} & \tfrac{b^2}{2ab-1} & \tfrac{b}{2ab-1}    \\
 a^2 & ab & a  \\
 a & b & 1
\end{pmatrix} ,
\end{equation}
with the corresponding loop group element
\begin{equation} \label{rk1}
g(\lambda) = I + \frac{2}{\lambda^3-\alpha^3}
\begin{pmatrix}
 \tfrac{\alpha^3 ab}{2ab-1} & \tfrac{\alpha \lambda^2 b^2}{2ab-1}   &  \tfrac{\alpha^2 \lambda b}{2ab-1}   \\
 \alpha^2 \lambda a^2 & \alpha^3 ab & \alpha \lambda^2 a  \\
 \alpha \lambda^2 a & \alpha^2 \lambda b &  \alpha^3
\end{pmatrix} ;
\end{equation}
and the \textbf{rank $2$ type} (\emph{thus the matrix \eqref{abc}
has rank $1$}) is as follows:
\begin{equation}
A= \frac{1}{3} \begin{pmatrix}
 \frac{ab-1}{2ab-1} & \frac{-b^2}{2ab-1} & \frac{b}{2ab-1}    \\
 a^2 & 1-ab & -a \\
 -a & b & 0
\end{pmatrix} ,
\end{equation}
with the corresponding loop group element
\begin{equation} \label{rk2}
g(\lambda) = I + \frac{2}{\lambda^3-\alpha^3}
\begin{pmatrix}
 \frac{\alpha^3(ab-1)}{2ab-1} & \frac{-\alpha \lambda^2 b^2}{2ab-1} & \frac{\alpha^2 \lambda b}{2ab-1}    \\
 \alpha^2 \lambda a^2 & \alpha^3 (1-ab) & -\alpha \lambda^2 a    \\
 -\alpha \lambda^2 a & \alpha^2 \lambda b & 0
\end{pmatrix}.
\end{equation}
Both types must meet the constraint: $2ab \neq 1$.
\begin{rem}\label{cpldet}
We compute that $\det g$ is
$[(\lambda^3+\alpha^3)/(\lambda^3-\alpha^3)]^{\mathrm{rank}(A)}$,
i.e., only depending on the poles and the rank of the residues. A
scaling by $ (\det g)^{-1/3} $ will make them lie in $\SL(3, \C)$,
though not rational any more.  Since the scaling does not affect the
factorization and the dressing action, we will ignore this step
henceforth.
\end{rem}

Let $l:=(a,b,1)$, $\ell$ be the line $\C \cdot l$, and introduce the
following `cone':
\[
\Delta:=\{\, (z_1,z_2,z_3) \in \C^3 ~ | ~ 2z_1 z_2=z_3^{2}, \textrm{
or }z_3=0 \, \}  .
\]
Then $2ab \neq 1$ is equivalent to $\ell \nsubseteq \Delta$. We
observe that $\ell^t = \Ima(\Res_{\alpha}g^t)$ for rank $1$ type and
$\ell^t = \Ker(\Res_{\alpha}g \, P)$ for rank $2$  type (here
$\ell^t$ means $\C \cdot l^t$). Conversely, a line $\ell$ not in
$\Delta$ determines the residue $A$ and thus the simple element $g$
uniquely in both types. Henceforth we always use $g_{\alpha, \ell}$
to denote the rank $1$ type element \eqref{rk1} and use $h_{\alpha,
\ell}$ to denote the rank $2$ type element \eqref{rk2}. We have
proved the following theorem:
\begin{thm}
The simplest rational element in $\Lambda_{-}^{\sigma}\GL(3, \C)$ is
either $g_{\alpha, \ell}$ of rank $1$ type \eqref{rk1} or
$h_{\alpha, \ell}$ of rank $2$ type \eqref{rk2},  where $\alpha \in
\C_{\ast}$ and $\ell \nsubseteq \Delta$.
\end{thm}

Imposing the $\tau$-reality condition $\overline{g(\bar{\lambda})} =
g(\lambda)$ on both types, we obtain that one pole, say $\alpha$,
has to be real and so is the residue $A$ there. It is convenient in
this case to let $\ell$ denote the real line $\R \cdot (a,b,1)$ in
$\R^3$ and let $\Delta_0$ denote $\Delta \cap \R^3$. We then have:
\begin{cor}
The simplest rational element in $\Lambda_{-}^{\tau,\sigma}\GL(3,
\C)$ is either $g_{\alpha, \ell}$ of rank $1$ type or $h_{\alpha,
\ell}$ of rank $2$ type,  where $\alpha \in \R_{\ast}$ and the real
line $\ell \nsubseteq \Delta_0$.
\end{cor}

We are ready to compute the dressing action of these simple
elements.

\begin{lem}\label{lemdress}
Let $g_{\alpha,\ell}\in \Lambda^{\tau, \sigma}_{-} \GL(3,\C)$ as in
\eqref{rk1}, and $F \in \Lambda^{\tau, \sigma}_{+} \SL(3,\C)$. If \
$\tilde{\ell} := \ell F(\alpha) \nsubseteq\Delta_0$, then \
$g_{\alpha,\ell} \cdot F$ \ can be factored uniquely as
\[
g_{\alpha,\ell} \cdot F = \tilde{F} \cdot g_{\alpha,\tilde{\ell}}
~\in~ \Lambda^{\tau, \sigma}_{+} \SL(3,\C) \times \Lambda^{\tau,
\sigma}_{-} \GL(3,\C) .
\]
\end{lem}
\begin{proof}
It suffices to prove that $\tilde{F} := g_{\alpha,\ell} \cdot F
\cdot g_{\alpha,\tilde{\ell}}^{-1}$ lies in $\Lambda^{\tau,
\sigma}_{+} \SL(3,\C)$. Since $\tilde{F}$ satisfies the reality
conditions \eqref{eqrea} and is holomorphic in $\C_{\ast}$ except
for possible simple poles coming from the poles of $g_{\alpha,\ell}$
and $g_{\alpha,\tilde{\ell}}^{-1}$\ , we only need to prove that the
residues of $\tilde{F}$ are zero at both $\alpha$ and $-\alpha$. But
\[
\mu (~ g_{\alpha,\ell}(\lambda) ~) = g_{\alpha,\ell}(- \lambda) \iff
P = g_{\alpha,\ell}(\lambda) ~ P ~ g_{\alpha,\ell}(-\lambda)^{t}
 ~ ,
\]
whose residue is zero at $\alpha$ implies that $\ell P
g_{\alpha,\ell}(-\alpha)^{t} = 0$, or equivalently
$g_{\alpha,\ell}(-\alpha) P \ell^{t} = 0$. These two equations are
also true for $\tilde{\ell}$. Therefore $(a,b,1)F(\alpha) \in
\tilde{\ell}$ and the special form of $A$ in \eqref{res1} imply that
\[
\Res_{\alpha}\tilde{F}=2\alpha~A~F(\alpha)~P~
g_{\alpha,\tilde{\ell}}(-\alpha)^{t}~P=0 ,
\]
and $F(-\alpha) P\, (\tilde{a}, \tilde{b}, 1)^{t}\in
[F(-\alpha)PF(\alpha)^{t}] \, \ell^t=P \ell^t$ implies that
\[
\Res_{-\alpha}\tilde{F}=-2\alpha~g_{\alpha,\ell}(-\alpha)~F(-\alpha)~P~\tilde{A}^{t}~P=0.
\]
The proof is completed once we notice that $\det\tilde{F}=1$ by
Remark \ref{cpldet}.
\end{proof}

\begin{thm}\label{main}
The dressing action of rank $1$ type $g_{\alpha,\ell}$ on the affine
frames $F(x,y,\lambda)$ of proper indefinite affine spheres gives
the classical Tzitz\'eica transformation, provided an open condition
that $\ell F(x,y,\alpha)\nsubseteq \Delta_0$. The dressing action of
rank $2$ type $h_{\alpha,\ell}$ gives the dual transformation.
\end{thm}
\begin{proof}
By Lemma \ref{lemdress}, for fixed $(x, y)$, we have the
factorization
\[
g_{\alpha,\ell}(\lambda) \cdot F(x,y, \lambda) = \tilde{F} (x, y,
\lambda) \cdot g_{ \alpha, \tilde{\ell} }(\lambda)
\]
with $\tilde{\ell} = \ell F(x,y,\alpha)$. From $F(x,y,
\alpha)=\left( (X_{\alpha})_x / \alpha, ~ \alpha (X_{\alpha})_y /h,
~ X_{\alpha}  \right)$, we get
\[
 (a,b,1) ~ F(x,y,\alpha) = ( ~ \phi_x / \alpha, ~ \alpha \phi_y /h, ~\phi~)\in \tilde{\ell},
\]
where $\phi:=(a,b,1) X_{\alpha}$ is a scalar solution of
\eqref{eqmap} with parameter $\alpha$. Note that a constant scaling
of $\phi$ does not change Tzitz\'eica  transformation \eqref{eqtzt},
and the solution space of the linear system \eqref{eqmap} has
dimension $3$. Therefore, by varying $\ell$, $\phi$ can be generic
scalar solution up to a constant multiple.

By the discussion at the beginning of this section, the third column
of $\tilde{F}$  produces new affine sphere, so does the affine
transformation of it by $g_{\alpha,\ell}^{-1}$\ :
\begin{eqnarray*}
\hat{X} &:=& g_{\alpha,\ell}^{-1}(\tilde{F})_{3} \\
&=& \left(F(\lambda) ~ g_{\alpha,\tilde{\ell}}(\lambda)^{-1} \right)_{3} \\
&=& \left(F(\lambda) ~ P ~ g_{\alpha,\tilde{\ell}}(- \lambda)^{t} ~ P \right)_{3}  \\
&\overset{\eqref{rk1}}{=}& \left(~ X_x / \lambda, ~ \lambda X_y / h,
~ X \right) ~ \cdot ~ \frac{1}{ \lambda^{3} + \alpha^{3} } ~ \cdot ~
\begin{pmatrix}
  2 \alpha^{3} \lambda   \phi_{y} / (h \phi)  \\
  - 2 \lambda^{2}  \phi_{x} / \phi   \\
 \lambda^{3} - \alpha^{3}
\end{pmatrix} \\
&=& \frac{\displaystyle  (\lambda^{3}  - \alpha^{3} ) h X - 2
\lambda^{3}  (\ln\phi)_{x} X_{y} + 2 \alpha^{3} (\ln\phi)_{y} X_{x}
}{\displaystyle  (\lambda^{3} + \alpha^{3} ) h }  ~ .
\end{eqnarray*}
The corresponding solution to Tzitz\'eica equation is given by
\[
\hat{h} = \hat{X}_{xy} / \hat{X} = h - 2 (\ln \phi)_{xy} .
\]
This is exactly the classical Tzitz\'eica transformation
\eqref{eqtzt} with
\[
\gamma = \lambda^{3}, \qquad \gamma_{1} = \alpha^{3} , \qquad
\phi_{1} = \phi .
\]

In rank $2$ type case, there is a similar factorization
$h_{\alpha,\ell} \cdot F = \tilde{F} \cdot h_{\alpha,\tilde{\ell}}$
when $\tilde{\ell}:=\Ker^t(AF(\alpha)P) \nsubseteq\Delta_0$. We omit
the details and present the corresponding transformation on affine
spheres:
\[
\tilde{X} = \frac{\displaystyle  (\lambda^{3}  + \alpha^{3} ) h X -
2 \lambda^{3}  (\ln\phi)_{x} X_{y} - 2 \alpha^{3} (\ln\phi)_{y}
X_{x} }{\displaystyle  (\lambda^{3} - \alpha^{3} ) h } ,
\]
where $\phi$ is the same scalar solution as the rank $1$ type case.
We see that
\[
\hat{X}_{\lambda} = - \tilde{X}_{-\lambda} ~ ,
\]
i.e., $- \tilde{X}$ gives the dual of $\hat{X}$. This completes the
proof.
\end{proof}

\medskip
\section{Permutability theorem and complex Tzitz\'eica transformations}

Let us briefly review the classical description of the permutability
theorem. In Theorem \ref{thmctt}, let $\phi_{1}$, $\phi_{2}$ be the
scalar solution of \eqref{eqmap} with parameter $\gamma_{1}$,
$\gamma_{2}$ respectively. Then using $\phi_{1}$ to apply
Tzitz\'eica transformation on $h$, we get a new solution to
Tzitz\'eica equation:
\[
h_{1}  := h - 2 (\ln \phi_{1})_{xy} ~.
\]
Applying Tzitz\'eica transformation \eqref{eqtzt}  to $(\phi_{2},
\gamma_{2})$, we obtain
\begin{equation}\label{phi12}
\phi_{12} := \frac{\displaystyle  (\gamma_{2} - \gamma_{1} ) h
\phi_{2} - 2 \gamma_{2} (\ln\phi_{1})_{x} (\phi_{2})_{y} + 2
\gamma_{1} (\ln \phi_{1})_{y} (\phi_{2})_{x} }{\displaystyle
(\gamma_{2} + \gamma_{1} ) h }
\end{equation}
as a scalar solution to \eqref{eqmap} with new $h_{1}$ and parameter
$\gamma_{2}$. Therefore we can use $\phi_{12}$ to apply Tzitz\'eica
transformation again on the new $h_{1}$, i.e.,
\[
h_{12}  = h_{1} - 2 (\ln \phi_{12})_{xy}
\]
will give another solution to Tzitz\'eica equation. In this two step
iteration, we may interchange the roles of $\phi_{1}$ and $\phi_{2}$
to obtain $h_{21}$ as another new solution. The permutability
theorem claims $h_{12} = h_{21}$. Similarly we can apply this two
step iteration to the affine sphere $X$ to obtain $X_{12}$ and
$X_{21}$ respectively, and the equality $X_{12}= X_{21}$ still
holds.

We will give a group point of view to this permutability theorem.

\begin{lem}\label{lemperm}
Let $g_{\alpha_{i},\ell_{i}}(\lambda)$ ($i=1,2$) be of rank $1$ type
with $\alpha_1^3 \neq \pm\alpha_2^3$. If both $\tilde{\ell}_1 :=
\ell_1  g_{\alpha_2,\ell_2}(\alpha_1)^{-1}$ and $\tilde{\ell}_2 :=
\ell_2  g_{\alpha_1,\ell_1}(\alpha_2)^{-1}$ are not in $\Delta$,
then
\begin{equation}\label{perm}
g_{\alpha_{2},\tilde{\ell}_{2}} ~ g_{\alpha_{1},\ell_{1}} =
g_{\alpha_{1},\tilde{\ell}_{1}} ~ g_{\alpha_{2},\ell_{2}} ~ .
\end{equation}
\end{lem}
\begin{proof}
It is equivalent to prove that $f:=g_{\alpha_{1},\tilde{\ell}_{1}}
g_{\alpha_{2},\ell_2} g_{\alpha_{1},\ell_{1}}^{-1}$ equals
$g_{\alpha_{2},\tilde{\ell}_{2}}$. First of all, they both are
rational elements in the group $\Lambda^{\sigma}_{-} \GL(3,\C)$. It
suffices to prove that their poles and residues are the same.

Let $l_i=(a_i,b_i,1)$ and $\tilde{l}_i=(\tilde{a}_i,\tilde{b}_i,1)$
span $\ell_i$ and $\tilde{\ell}_i$ respectively. Similar to the
proof of Lemma \ref{lemdress}, we compute that
\[
\Res_{\alpha_1}f = 2\alpha_1
~\widetilde{A}_1~g_{\alpha_2,\ell_2}(\alpha_1)~P~
g_{\alpha_1,\ell_1}(-\alpha_1)^{t}~P =0
\]
since $\tilde{l}_1 g_{\alpha_2,\ell_2}(\alpha_1) \in \ell_1$, and
\[
\Res_{-\alpha_1} f =-2\alpha_1 ~
g_{\alpha_{1},\tilde{\ell}_{1}}(-\alpha_1) ~
g_{\alpha_2,\ell_2}(-\alpha_1)~P~A_1^t~P =0
\]
since $g_{\alpha_2,\ell_2}(-\alpha_1) P\, l_1^t \in [
g_{\alpha_2,\ell_2}(-\alpha_1)\, P\, g_{\alpha_2,\ell_2}(\alpha_1)^t
] \, \tilde{\ell}^t_1 = P \, \tilde{\ell}_1^t$. Thus $f$ has only
three simple poles $\{ \alpha_2, \epsilon^{2} \alpha_2, \epsilon^{4}
\alpha_2 \}$, same as $g_{\alpha_2,\tilde{\ell}_{2}}$.

Now due to the $\nu$ reality condition in \eqref{eqrea}, we only
need to prove that their residues at $\alpha_2$ are the same. But
\begin{eqnarray*}
\Ima \left(\Res_{\alpha_2}f^t \right) &=& \Ima \left(\left(g_{\alpha_1,\ell_1}(\alpha_2)^{-1}\right)^t A_2^t ~ g_{\alpha_1,\tilde{\ell}_1}(\alpha_2)^t \right) \\
&=&  \left(g_{\alpha_1,\ell_1}(\alpha_2)^{-1}\right)^t \ell_2^t \\
&=& \tilde{\ell}_2^t \quad = \quad \Ima \left(\Res_{\alpha_2}
g_{\alpha_2,\tilde{\ell}_2}^t \right) .
\end{eqnarray*}
Then $\Res_{\alpha_2}f$ must be the same as $\Res_{\alpha_2}
(g_{\alpha_2,\tilde{\ell}_2})$ since rank $1$ type residue is
uniquely determined by the above image. This  completes the proof.
\end{proof}

\begin{eg}\label{eg0}
Choose two nonzero poles $\alpha_{1}, \alpha_{2}$ such that
$\alpha_1^3\neq \pm\alpha_2^3$, and let $\ell_i = \C \cdot
(0,b_i,1)$ for $i=1,2$. Then $g_{\alpha_{2},\tilde{\ell}_{2}} ~
g_{\alpha_{1},\ell_{1}} = g_{\alpha_{1},\tilde{\ell}_{1}} ~
g_{\alpha_{2},\ell_{2}}$ holds for $\tilde{\ell}_i = \C \cdot (0,
\tilde{b}_i, 1)$ with
\[
\tilde{b}_1 = \frac{ (\alpha_1^3 + \alpha_2^3) b_1 - 2 \alpha_1
\alpha_2^2 b_2 }{ (\alpha_1^3 - \alpha_2^3) } , \qquad \tilde{b}_2 =
\frac{ 2 \alpha_1^2 \alpha_2 b_1 - (\alpha_1^3 + \alpha_2^3) b_2 }{
(\alpha_1^3 - \alpha_2^3) } .
\]
\end{eg}

\smallskip

\begin{thm}
Use the same notation and the factorization formula \eqref{perm} in
Lemma \ref{lemperm}. Let all $\alpha_i,\ell_i$ be real. Let
$F_1:=g_{\alpha_1,\ell_1} \ast F$ and $F_2:=g_{\alpha_2,\ell_2} \ast
F$, where $\ast$ is the dressing action on the frames
$F(x,y,\lambda)$ of affine spheres. Then the following holds and
implies the classical permutability theorem:
\begin{equation}\label{gpperm}
\begin{split}
  & ( g_{\alpha_{2},\tilde{\ell}_{2}} g_{\alpha_{1},\ell_{1}} ) \ast F =
g_{\alpha_{2},\tilde{\ell}_{2}} \ast F_{1} \quad (~ =: F_{12} ~)  \\
= ~ & ( g_{\alpha_{1},\tilde{\ell}_{1}} g_{\alpha_{2},\ell_{2}} )
\ast F = g_{\alpha_{1},\tilde{\ell}_{1}} \ast F_{2} \quad (~ =:
F_{21} ~) .
\end{split}
\end{equation}
\end{thm}
\begin{proof}
Because the dressing action is a group action, \eqref{gpperm}
certainly holds by \eqref{perm}. Let $l_i=(a_i,b_i,1)$ and
$\tilde{l}_i=(\tilde{a}_i,\tilde{b}_i,1)$ span $\ell_i$ and
$\tilde{\ell}_i$ respectively. From the proof of Theorem \ref{main},
$F_1:=g_{\alpha_1,\ell_1} \ast F$ means in classical terms the
following: Tzitz\'eica transformation via $\phi_{1} := l_1
(F(\alpha_1))_{3} $ on $X=(F)_{3}$ gives a new affine sphere $X_{1}
:= g_{\alpha_{1},\ell_{1}}^{-1} (F_{1})_{3}$. Therefore $F_{12} :=
g_{\alpha_{2},\tilde{\ell}_{2}} \ast F_{1} $ implies that
Tzitz\'eica transformation via $\phi_{12} := \tilde{l}_2
(F_1(\alpha_2))_{3} $ on $(F_1)_{3} = g_{\alpha_{1},\ell_{1}} X_{1}$
gives a new affine sphere $ g_{\alpha_{2},\tilde{\ell}_{2}}^{-1}
\left(F_{12} \right)_{3}$. We observe that
\[
\phi_{12} :=\tilde{l}_2 (F_1(\alpha_2))_{3}  = c_0 l_2 \,
g_{\alpha_1,\ell_1}(\alpha_2)^{-1} (F_1(\alpha_2))_{3} =
 c_{0}  l_2 \, X_1(\alpha_2)  ,
\]
which coincides with the classical formula \eqref{phi12} except for
a negligible constant $c_{0}$ when we plug in $\phi_{2} := l_{2}\,
X(\alpha_2)$.

So Tzitz\'eica transformation via $\phi_{12}$ on $X_{1} =
g_{\alpha_{1},\ell_{1}}^{-1} (F_1)_3$ produces
\[
X_{12} = g_{\alpha_{1},\ell_{1}}^{-1} \left(
g_{\alpha_{2},\tilde{\ell}_{2}}^{-1} F_{12} \right)_{3}= \left(
g_{\alpha_{1},\ell_{1}}^{-1} g_{\alpha_{2},\tilde{\ell}_{2}}^{-1}
F_{12} \right)_{3} ~ .
\]
Similarly $X_{21}=\left( g_{\alpha_{2},\ell_{2}}^{-1}
g_{\alpha_{1},\tilde{\ell}_{1}}^{-1} F_{21} \right)_{3}$. Therefore,
by \eqref{perm} and \eqref{gpperm}, we obtain the classical
permutability theorem: $X_{12} = X_{21}$, which automatically
implies $h_{12} = h_{21}$ for the corresponding affine metrics.
\end{proof}

There is some rational element in
$\Lambda^{\tau,\sigma}_{-}\GL(3,\C)$ which has $6$ simple poles but
none of them are real. The poles must form two conjugate triples as
$\{ \alpha, \epsilon^{2} \alpha, \epsilon^{4} \alpha \}$ and
$\{\bar{\alpha}, \epsilon^{2} \bar{\alpha}, \epsilon^{4}
\bar{\alpha} \}$, where we may assume $0<\arg(\alpha)<\pi/3$ without
loss of generality. So such element is not product of \emph{real}
rank $1$ or $2$ type  elements. In fact, we can use Lemma
\ref{lemperm} to construct them:
\begin{prop}\label{cplpole}
Let $\alpha \in \C_{\ast} $ with $\arg(\alpha) \in (0,\pi/6) \cup
(\pi/6,\pi/3)$. Let $\ell \nsubseteq \Delta $. If $\ell^{\ast}
:=\bar{\ell} \cdot g_{\alpha,\ell}(\bar{\alpha})^{-1} \nsubseteq
\Delta $, then $f_{\alpha,\ell} := g_{\bar{\alpha}, \ell^{\ast} } ~
g_{\alpha, \ell} \in \Lambda^{\tau,\sigma}_{-}\GL(3,\C) $.
\end{prop}
\begin{proof}
We first observe that $f_{\alpha,\ell}$ lies in
$\Lambda^{\sigma}_{-}\GL(3,\C)$. It suffices to verify
$\overline{f(\bar{\lambda})} = f(\lambda)$, which is
\[
g_{\alpha, \overline{\ell^{\ast}} } ~~ g_{\bar{\alpha}, \bar{\ell} }
~=~ g_{\bar{\alpha}, \ell^{\ast}} ~~   g_{\alpha, \ell } .
\]
Since $\ell^{\ast} = \bar{\ell} \cdot
g_{\alpha,\ell}(\bar{\alpha})^{-1}\nsubseteq \Delta$ implies
$\overline{\ell^{\ast}} = \ell \cdot
g_{\bar{\alpha},\bar{\ell}}(\alpha)^{-1} \nsubseteq \Delta$ and
$\arg(\alpha)\neq \pi/6$ implies $\alpha^3\neq -\bar{\alpha}^3$, the
above factorization holds by Lemma \ref{lemperm}.
\end{proof}

The dressing action of $f_{\alpha,\ell}$ on affine spheres can be
viewed as the composition of two `conjugate' complex Tzitz\'eica
transformations, which produces a real solution in the end. The
Permutability Theorem \ref{gpperm} can be applied to compute this
action. Solutions from this construction are often called breather
type  solutions.

It is not hard to show, by a similar residue calculus as before,
that any rational element with $6$ simple poles as above can be
constructed from Proposition \ref{cplpole}. What  would be much
messier, if not harder, to prove is that the subgroup of all
rational elements in $\Lambda^{\tau,\sigma}_{-}\GL(3,\C)$ is
generated by $g_{\alpha,\ell}$'s,  $f_{\alpha,\ell}$'s, and their
rank $2$ type brothers. This subgroup can then be regarded as the
group of Tzitz\'eica transformations on affine spheres. We will
leave this interesting problem for future study.

\medskip
\section{Basic examples}

In this section, we use $x, y, z$ as the standard $\R^{3}$
coordinates to represent the immersion $X$, and use $u, v$ to denote
the asymptotic coordinates of the affine spheres.

\begin{eg}[The \textbf{vacuum} solution]
The vacuum solution to Tzitz\'eica equation  (see also
\cite{DorEit01}, \cite{RogSch94}) is $\omega_{0} \equiv 0$ (or
$h_{0} \equiv 1$). One can integrate \eqref{eqframe} to obtain the
whole family of frames. The Cartesian equation of the surface is
then obtained by the determinant:
\[
x^3+y^3+z^3-3xyz=1.
\]
Note that it is independent of the parameter $\lambda$. So this
family is really a family of parametrizations of the same affine
sphere.

A general scalar solution of system \eqref{eqmap} with parameter
$\gamma =  \lambda^{3} $ is:
\begin{equation}\label{eqsca}
\phi(\lambda) = c_{0} R(\lambda) + c_{1} R(\epsilon^{2} \lambda)+
c_{2} R(\epsilon^{4} \lambda),
\end{equation}
where $R(\lambda):=\exp(\lambda u+\lambda^{-1}v)$.

Therefore we may choose the following asymptotic parametrizations of
the vacuum affine sphere after certain affine transformation:
\[
X_{0}(u,v,\lambda) =
\begin{pmatrix}
 \exp[-(\lambda u + \lambda^{-1} v )/2] ~\cos[\sqrt{3} (\lambda u - \lambda^{-1} v )/2 ] \\
 \exp[-(\lambda u + \lambda^{-1} v )/2] ~\sin[\sqrt{3} (\lambda u - \lambda^{-1} v )/2 ] \\
  \frac{2}{3\sqrt{3}} ~ \exp(\lambda u + \lambda^{-1} v )
\end{pmatrix},
\]
which is a surface of revolution (Note that Jonas has classified all
affine spheres of revolution in \cite{Jon21} using elliptic
functions).
\end{eg}

\begin{eg}[The \textbf{one-soliton} solution]
Apply Tzitz\'eica transformation to the vacuum solution we obtain
the one-soliton solution $h_{1}$. By \eqref{eqsca},
$\phi_{1}=\phi(\lambda_{1})$ is a scalar solution of system
\eqref{eqmap} with parameter $\gamma_{1} =  \lambda_{1}^{3} $. It is
real when $c_{0}\in\R$, $c_{1} = \overline{c_2}$, and $\lambda_{1}
\in \R_{\ast}$. Compute $h_{1}=1-2(\ln \phi_{1})_{uv}$:
\[
h_{1} = 1 - \frac{\displaystyle 6\beta_{0} \exp(3s_1/2) \cos \left(
\sqrt{3} \ t_1 / 2 + \theta_{0} \right) + 1.5 }{\displaystyle \left[
\beta_{0} \exp(3s_1/2) + \cos \left( \sqrt{3} \ t_1 / 2 + \theta_{0}
\right)  \right]^{2} },
\]
where $c_1=\rho_{0} e^{ {\rm i}\theta_{0} } $, $\beta_{0} = c_{0} /
(2 \rho_{0})$, $s_1=\lambda_{1} u + \lambda_{1}^{-1}v $, and
$t_1=\lambda_{1} u - \lambda_{1}^{-1}v $. The family of affine
apheres $X_{1}(u,v,\lambda)$ has a long expression given by
\eqref{eqtzt}.

When $\beta_{0} =0$ (i.e. $c_{0}=0$), we have the special solution:
\begin{equation}\label{oneh}
h_{1} = 1 - 1.5 \sec^{2} \left[ \sqrt{3} \ (\lambda_{1} u -
\lambda_{1}^{-1}v) / 2 + \theta_{0}  \right]   .
\end{equation}
We give explicit formula for this family of affine spheres:
\begin{equation}\label{kelch}
\begin{split}
X_{1}(u,v,\lambda) & =
\frac{(\lambda^3-\lambda_1^3)}{(\lambda^3+\lambda_1^3)}X_{0}(u,v,\lambda)
+ \frac{ \sqrt{3} \lambda \lambda_1 \tan\left( \sqrt{3} t_1 /2 +
\theta_0  \right)}{(\lambda^3+\lambda_1^3)} \cdot  \\ {} & {}
\begin{pmatrix}
 e^{-\frac{s}{2}}\left[ \lambda \cos\left( \sqrt{3} t /2 + 4\pi /3  \right) + \lambda_1 \cos\left( \sqrt{3} t /2 + 2\pi /3 \right) \right] \\
 e^{-\frac{s}{2}}\left[ \lambda  \sin\left( \sqrt{3} t /2 + 4\pi /3  \right) + \lambda_1 \, \sin\left( \sqrt{3} t /2 + 2\pi /3  \right) \right] \\
  2 e^{s} (\lambda+\lambda_1) / (3\sqrt{3})
\end{pmatrix} .
\end{split}
\end{equation}

Note that $\phi_1$ need not be real to produce real $h_1$. For
example, \eqref{oneh} will be a real hyperbolic function solution
when $\lambda_{1}$ and $\theta_{0}$ are pure imaginary. In this case
the real (or imaginary) part of \eqref{kelch} still produces affine
spheres, among which are the stationary and traveling one-soliton
affine sphere shown in \cite{Sch01}. Some pictures have already been
shown in \cite{RogSch94} and \cite{Sch01}. In terms of dressing
action, $\lambda_1$ is the pole of some rank $1$ type simple element
$g_{\lambda_1,\ell}$. So dressing actions of $g_{\lambda_1,\ell}$
with a pure imaginary pole (or a pole whose argument is $\pm \pi/6$)
may also produce new real affine spheres sometime.
\end{eg}

\medskip
\section{Acknowledgments}
Part of this paper was in the author's thesis at Northeastern
University and he would like to express the deepest gratitude to his
advisor Chuu-Lian Terng, for her help and encouragement.  The author
would like to thank Franz Pedit for pointing out an error about the
associated family of affine spheres, when the author presented a
draft of this paper in the $2005$ AMS meeting at Lubbock, Texas.
Thanks also go to the organizers of this meeting Josef F.
Dorfmeister and Hongyou Wu for the invitation and helpful
discussions.

The research was supported in part by Postdoctoral Fellowship of the
Mathematical Sciences Research Institute.


\bibliographystyle{alpha}

\end{document}